\documentclass[a4paper, 12pt]{article}

\usepackage{amsfonts, amsmath, amsthm}
\usepackage{amssymb}
\usepackage{latexsym}
\usepackage[dvips]{graphicx}
\usepackage{color}

\theoremstyle{plain}
\newtheorem{thm}{Theorem}[section]
\newtheorem{prop}{Proposition}[section]
\newtheorem{lem}{Lemma}[section]
\newtheorem{cor}{Corollary}[section]
\theoremstyle{definition}

\theoremstyle{definition}
\newtheorem{rmk}{Remark}[section]

\numberwithin{equation}{section}

\newcommand{\Z}{\mathbb{Z}}

\newcommand{\R}{\mathbb{R}}

\newcommand{\Sph}{\mathbb{S}}

\newcommand{\op}[1]{\mathcal{#1}}

\newcommand{\pa}{\partial}
\newcommand{\eps}{\varepsilon}
\newcommand{\jb}[1]{\langle #1 \rangle}

\newcommand{\LR}[1]{\left| #1 \right|}
\newcommand{\lr}[1]{\left( #1 \right)}

\begin{document}
\title{
Non-decay of the energy for a system of semilinear wave equations
}
\author{
          Yoshinori Nishii\thanks{
              Department of Mathematics, Graduate School of Science,
              Osaka University.
              1-1 Machikaneyama-cho, Toyonaka,
              Osaka 560-0043, Japan.
              (E-mail: {\tt y-nishii@cr.math.sci.osaka-u.ac.jp})
             }
}
\date{\today }
\maketitle


\noindent{\bf Abstract:}\ 
We consider the global Cauchy problem for a two-component system 
of cubic semilinear wave equations in two space dimensions.
We give a criterion for large time non-decay of the energy 
for small amplitude solutions in terms of 
the radiation fields associated with the initial data.
\\

\\



\section{Introduction}  \label{sec_intro}
This paper is intended to be a continuation of \cite{NS}. 
We are interested in large time behavior of solutions to the Cauchy problem for
\begin{align}
\left\{\begin{array}{ll}
  \begin{array}{l}
\Box u_1 = - (\pa_t u_2)^2 \pa_t u_1, \\
\Box u_2 = - (\pa_t u_1)^2 \pa_t u_2,
  \end{array}
  & (t,x)\in (0,\infty)\times\R^2,
 \end{array}\right.
\label{eq_2}
\end{align}
with the initial condition
\begin{align}
u_j(0,x)=\eps f_j(x), \ \pa_t u_j(0,x)=\eps g_j(x), \qquad x\in\R^2,\ j=1,2,
\label{data}
\end{align}
where 
$\Box =\pa_t^2-\triangle=\pa_0^2-(\pa_1^2+\pa_2^2)$ with
$\pa_0=\pa_t=\pa/\pa t$, 
$\pa_1=\pa/\pa x_1$,
$\pa_2=\pa/\pa x_2$.
We assume that $f=(f_1,f_2), g=(g_1,g_2)$ are compactly supported $C^{\infty}$ functions
on $\R^2$ and $\eps$ is a small positive parameter.
The following result has been obtained in \cite{NS}:
\begin{prop}[\cite{NS}] \label{prop_NS1}
The Cauchy problem \eqref{eq_2}--\eqref{data} admits a unique 
global $C^{\infty}$ solution if $\eps$ is suitably small. 
Moreover, there exists $(f^+,g^+)\in \dot H^1(\R^2)\times L^2(\R^2)$ such that
\begin{align} \label{asympt_free}
\lim_{t\to+\infty}\| u(t)-u^+(t) \|_E=0,
\end{align}
where 
$u^+=(u^+_1,u^+_2)$ solves the free wave equation 
$\Box u^+=0$ with $(u^+,\pa_t u^+)|_{t=0}=(f^+,g^+)$ 
and the energy norm $\|\cdot\|_E$ is defined by
\[
\| \phi(t) \|_E^2:=\frac{1}{2}\int_{\R^2}|\pa \phi(t,x)|^2 \, dx.
\]
\end{prop}


As emphasized in \cite{NS}, this is not a trivial result because 
cubic nonlinearities must be regarded as long-range perturbation 
for two-dimensional wave equations in general.
Note that the global existence part of this assertion follows from the
earlier result \cite{KMatsS} directly 
and that the system \eqref{eq_2} possesses two conservation laws
\[
\frac{d}{dt}\lr{\|u_1(t)\|_E^2+\|u_2(t)\|_E^2}
=-2\int_{\R^2}(\pa_tu_1(t,x))^2(\pa_tu_2(t,x))^2 \, dx
\]
and
\begin{align} \label{1-2}
\frac{d}{dt}\lr{\|u_1(t)\|_E^2-\|u_2(t)\|_E^2}=0.
\end{align}
However, these are not enough to conclude that the solution $u(t)$ is 
asymptotically free in the sense of \eqref{asympt_free}.
The novelty of the previous work \cite{NS} is to address this point, 
but we must say that the problem on (non-)triviality of $u^+$ is still obscure.
What we are going to do in the present work is to investigate it in more detail.
To make our concern clearer, let us focus on the special case $f_1=f_2$, $g_1=g_2$ first.
By the uniqueness of the solution, the problem is reduced 
to the single equation $\Box v = -(\pa_t v)^3$, 
whence we can adapt the result of \cite{KMuS}, \cite{KMatsS}, \cite{NST} to see that 
the total energy $\| u(t) \|_E$ decays like $O((\log t)^{-1/4+\delta})$ 
as $t\to+\infty$ for arbitrarily small $\delta>0$. In other words, 
\eqref{asympt_free} holds with the trivial free solution $u^+\equiv0$ in this case.
However, we remark that this is an exceptional case. 
Indeed, it follows from the conservation law \eqref{1-2} that at least one component 
$u_1$ or $u_2$ tends to a non-trivial free solution if $\|u_1(0)\|_E\neq\|u_2(0)\|_E$.
In other words, at least $u^+_1$ or $u^+_2$ given in Proposition~\ref{prop_NS1} does not 
vanish for generic initial data.
Moreover,  it is far from obvious whether both $u^+_1$ and $u^+_2$ can behave 
like non-trivial free solutions as $t\to+\infty$ in a certain case. 
This is why the problem on (non-)triviality of scattering state for \eqref{eq_2} 
is of our interest.

Before going further, let us review the strategy of the proof of Proposition~\ref{prop_NS1} briefly.
The key in \cite{NS} is to focus on the function $V=(V_1,V_2)$ defined by 
$V_j(t;\sigma,\omega)=U_j(t,(t+\sigma)\omega)$, where
$U_j(t,x)=\op{D}(|x|^{1/2}u_j(t,x))$, 
$\op{D}=2^{-1}(\pa_r-\pa_t)$, 
$\pa_{r}=\omega_{1}\pa_{1}+\omega_{2}\pa_{2}$, $r=|x|$ 
and $\omega =(\omega_1,\omega_2)=x/|x|$.
Roughly speaking, what have been seen in \cite{NS} is that 
the leading part of $\pa u_j(t,x)$ as $t\to+\infty$ can be given by
$|x|^{-1/2}\hat\omega(x) V_j(t;|x|-t,x/|x|)$, 
where $\pa u=(\pa_0 u,\pa_1 u,\pa_2 u)$ and 
$\hat \omega(x)=(\omega_0,\omega_1,\omega_2)=(-1,\omega)$.
Moreover, the evolution of $V=(V_1,V_2)$ can be characterized by the system
\begin{align} \label{profile}
\pa_{t} V_1
=\dfrac{-1}{2t} V_1V_2^2+ K_1, \qquad
\pa_{t} V_2
=\dfrac{-1}{2t} V_1^2V_2+ K_2,
\end{align}
where $K_1$ and $K_2$ are harmless remainder terms 
(see Proposition~\ref{prop_NS3} below).
These allow us to reduce the problem to investigating the behavior of 
the solution $V$ to \eqref{profile} as $t\to+\infty$.
As for the asymptotic behavior of $V$, we already know the following two propositions.
\begin{prop}[\cite{NS}] \label{prop_NS2}
The limit $V_j^+(\sigma,\omega)=\lim_{t\to+\infty}V_j(t;\sigma,\omega)$ exists for each 
fixed $(\sigma,\omega)\in\R\times\Sph^1$ and $j=1,2$, where $\Sph^1$ denotes 
the unit circle in $\R^2$.
Moreover, we have
\begin{align}\label{NS2}
\lim_{t\to\infty}\| \pa u_j(t,\cdot) - \hat\omega(\cdot)V_j^{+,\#}(t,\cdot) \|_{L^2(\R^2)}=0,
\end{align}
where 
$V_j^{+,\#}(t,x)=|x|^{-1/2}V^{+}_j(|x|-t,x/|x|)$.
\end{prop}
\begin{prop}[\cite{NS}] \label{prop_NS3}
Let $V^+_j$ be as above.
There exists a function $m:\R\times\Sph^1\to\R$ such that the following holds
for each $(\sigma,\omega)\in\R\times\Sph^1$:
\begin{itemize}
\item $m(\sigma,\omega)>0$ implies 
$V_1^+(\sigma,\omega)\neq0$ and $V_2^+(\sigma,\omega)=0$;
\item $m(\sigma,\omega)<0$ implies 
$V_1^+(\sigma,\omega)=0$ and $V_2^+(\sigma,\omega)\neq0$;
\item $m(\sigma,\omega)=0$ implies 
$V_1^+(\sigma,\omega)=V_2^+(\sigma,\omega)=0$.
\end{itemize}
We also have following expression of $m(\sigma,\omega)$:
\begin{align}\label{m}
m(\sigma,\omega)
=(V_1(t_{0,\sigma};\sigma,\omega))^2-(V_2(t_{0,\sigma};\sigma,\omega))^2
+ 2\int_{t_{0,\sigma}}^{\infty} \rho(\tau;\sigma,\omega) \, d\tau  ,
\end{align}
where 
$
t_{0,\sigma}=\max\{2,-2\sigma\}, 
$
\[
\rho(t;\sigma,\omega)
= V_1(t;\sigma,\omega)K_1(t;\sigma,\omega) 
- V_2(t;\sigma,\omega)K_2(t;\sigma,\omega), 
\]
\[
K_j(t;\sigma,\omega)=H_j(t,(t+\sigma)\omega), \qquad j=1,2, 
\]
\[
H_1(t,x)=\dfrac{1}{2}
\left( |x|^{1/2}(\pa_t u_2)^2(\pa_t u_1)+\dfrac{1}{t}U_2^2U_1 \right) 
-\dfrac{1}{8|x|^{3/2}}(4(x_1 \pa_2 - x_2 \pa _1)^2+1)u_1, 
\]
\[
H_2(t,x)=\dfrac{1}{2}
\left( |x|^{1/2}(\pa_t u_1)^2(\pa_t u_2)+\dfrac{1}{t}U_1^2U_2 \right) 
-\dfrac{1}{8|x|^{3/2}}(4(x_1 \pa_2 - x_2 \pa _1)^2+1)u_2.
\]
\end{prop}

From Proposition~\ref{prop_NS2}, we see that the vanishing of $V_j^+(\sigma,\omega)$
implies the energy decay of $u_j(t)$ (or, equivalently, the triviality of $u_j^+$).
Note also that the function $m(\sigma,\omega)$ appearing in Proposition~\ref{prop_NS3}
gives us much information on the vanishing of $V_j^+(\sigma,\omega)$.
Accordingly, it is natural to expect that the better understanding of $m(\sigma,\omega)$ 
provides us with more precise information on the energy decay. 
The aim of this paper is to specify the leading term of $m(\sigma,\omega)$ 
for sufficiently small $\eps$. 
As a consequence, we will find a criterion for non-decay of the energy 
for small amplitude solution to \eqref{eq_2}--\eqref{data} in the terms of the 
radiation fields associated with the initial data.
In particular, we will see that both $u_1(t)$ and $u_2(t)$ can behave like
non-zero free solutions as $t\to+\infty$ with a suitable choice of $(f,g)$.

\section{Main results}  \label{sec_main}
Before stating the main result, we introduce several notations.
For $\phi,\psi\in C_{0}^{\infty}$, we define 
the Friedlander radiation field $\mathcal{F}_0[\phi,\psi]$ by
\[
\op{F}_0[\phi,\psi](\sigma,\omega)
:=-\pa_{\sigma}\op{R}_2[\phi](\sigma,\omega)+\op{R}_2[\psi](\sigma,\omega),
\qquad (\sigma,\omega)\in\R\times\Sph^1,
\]
where
\[
\op{R}_2[\phi](\sigma,\omega)
:=\frac{1}{2\sqrt{2}\pi}\int_{\sigma}^{\infty} 
\frac{\op{R}[\phi](s,\omega)}{\sqrt{s-\sigma}}\, ds, \quad
\op{R}[\phi](s,\omega):=\int_{y\cdot\omega=s} \phi(y)\, dS_{y}.
\]
For simplicity, we write 
$\op{F}_j(\sigma,\omega)=\mathcal{F}_0[f_j,g_j](\sigma,\omega)$,  $j=1,2$.

Our main theorem is as follows.
\begin{thm}  \label{thm_main}
Let $m=m(\sigma,\omega)$ be the function given in Proposition~\ref{prop_NS3} 
and $0<\mu<1/10$. 
Then we have 
\begin{align}
m(\sigma, \omega) = \eps^2\left( \left(
\pa_{\sigma}\op{F}_1(\sigma,\omega)\right)^2 
- \left( \pa_{\sigma}\op{F}_2(\sigma,\omega) \right)^2 \right)
+O(\eps^{5/2-2\mu})
\end{align}
uniformly in $(\sigma,\omega)\in\R\times\Sph^1$ as $\eps \to +0$.
\end{thm}
This asymptotic expression of $m(\sigma,\omega)$ yields the following criterion for 
non-decay of $\| u_1(t)\|_E$ and $\| u_2(t)\|_E$.
\begin{cor}  \label{cor_1}
Suppose that there exist $(\sigma^{*},\omega^{*})$, $(\sigma_{*},\omega_{*})\in\R\times\Sph^1$ 
satisfying 
\begin{align}\label{1>2}
|\pa_{\sigma}\op{F}_1(\sigma^{*},\omega^{*})| 
> |\pa_{\sigma}\op{F}_2(\sigma^{*},\omega^{*})|
\end{align}
and
\begin{align*}
|\pa_{\sigma}\op{F}_1(\sigma_{*},\omega_{*})| 
< |\pa_{\sigma}\op{F}_2(\sigma_{*},\omega_{*})|,
\end{align*}
respectively. Then we have
$\lim_{t\to+\infty}\| u_1(t) \|_{E}>0$ and
$\lim_{t\to+\infty}\| u_2(t) \|_{E}>0$
for suitably small $\eps$.
\end{cor}

\begin{rmk}
From Corollary~\ref{cor_1}, we can construct the solution $u=(u_1,u_2)$ to 
\eqref{eq_2}--\eqref{data} whose energy does not decay.
Consequently, Corollary~\ref{cor_1} shows that both $u_1(t)$ and $u_2(t)$ 
can behave like non-trivial free solutions as $t\to+\infty$ in a certain case.
\end{rmk}

\begin{rmk}
In \cite{LNSS2}, analogous results have been obtained 
for the nonlinear Schr\"odinger system of the form
\begin{align*}
\left\{\begin{array}{ll}
  \begin{array}{l}
(i\pa_t +\frac{1}{2}\pa_x^2) u_1 = -i|u_2|^2  u_1, \\[5pt]
(i\pa_t +\frac{1}{2}\pa_x^2) u_2 = -i|u_1|^2  u_2,
  \end{array}
  & t>0,\ x\in\R.
 \end{array}\right.
\end{align*}
Many parts of our proof below are similar to those of \cite{LNSS2},
but we need several modifications to fit for the wave equation case.
\end{rmk}

The rest part of this paper is organized as follows.
In Section~\ref{sec_pre}, we introduce preliminary estimates for the 
small amplitude solution $u$ to \eqref{eq_2}--\eqref{data}.
We also recall basic estimates for solutions to the free wave equation
and the radiation fields associated with its initial data.
We prove Theorem~\ref{thm_main} and Corollary~\ref{cor_1} in Section~\ref{sec_proofs}.

\section{Preliminaries}  \label{sec_pre}

In this section, we collect various estimates which will be used 
in the next section.
First of all, we introduce several notations.
For $z\in \R$, we write $\jb{z} =\sqrt{1+|z|^{2}}$.
We define
\begin{align*}
 S := t \pa_{t} + x_{1}\pa_{1} + x_{2}\pa_{2}, \ 
 L_{1} := t \pa_{1} + x_{1} \pa_{t}, \ 
 L_{2} := t \pa_{2} + x_{2} \pa_{t}, \ 
 \Omega := x_1 \pa_2 - x_2 \pa _1,
\end{align*}
and we set
$\Gamma =(\Gamma_{j})_{0 \le j \le 6}
=(S , L_{1} , L_{2} , \Omega , \pa_{0} , \pa_{1} , \pa_{2})$.
For a multi-index 
$\alpha = (\alpha _{0},\alpha _{1},\cdots,\alpha _{6}) \in \Z^{7}_{+}$, 
we write 
$|\alpha |=\alpha _{0}+\alpha _{1}+\cdots+\alpha _{6}$ 
and 
$\Gamma^{\alpha}
=\Gamma_{0}^{\alpha_{0}}\Gamma_{1}^{\alpha_{1}}\cdots \Gamma_{6}^{\alpha_{6}}$, 
where $\Z_{+}=\{n\in\Z ; n\ge 0\}$.
We also define $|\, \cdot\,|_s$ and $\|\, \cdot\,\|_s$ by 
\begin{align*}
 |\phi (t,x)|_{s}=\sum _{|\alpha |\le s}|\Gamma^{\alpha }\phi (t,x)|, \qquad 
\| \phi (t,\cdot )\| _{s}=\sum _{|\alpha |\le s}\| \Gamma^{\alpha}\phi (t,\cdot) \| _{L^{2}(\R ^{2})}.
\end{align*}
From the argument of Section 3 in \cite{KMatsS}, we already know that the following 
estimates are satisfied by
the global small amplitude solution $u$ to \eqref{eq_2}--\eqref{data}.
\begin{lem}\label{lem_apriori}
Let $k\ge 4$, $0<\mu<1/10$ and $0<(8k+7)\nu<\mu$.
If $\eps>0$ is suitably small, then the solution 
$u$ to \eqref{eq_2}--\eqref{data} satisfies 
\begin{align}
&|u(t,x)|_{k+1} \le C\eps \jb{t+|x|}^{-1/2+\mu}, \label{apriori_u}\\
&|\pa u(t,x)| \le C\eps \jb{t+|x|}^{-1/2}\jb{t-|x|}^{\mu-1}, \label{apriori_pa_u}\\
&|\pa u(t,x)|_{k} \le C\eps \jb{t+|x|}^{-1/2+\nu}\jb{t-|x|}^{\mu-1}, \label{apriori_pa_u_k}
\end{align}
for $(t,x)\in [0,\infty)\times \R^2$ and
\begin{align}
\|\pa u(t) \|_{k}\le C\eps(1+t)^{\mu-\nu} \label{apriori_energy}
\end{align}
for $t\ge 0$, where $C$ is a positive constant independent of $\eps$.
\end{lem}

Next, we recall some estimates relevant to the free wave equation
\begin{align}
\left\{
  \begin{array}{ll}
\Box \phi = 0, \qquad (t,x)\in (0,\infty)\times\R^2,\\
  \begin{array}{l}
\phi(0)=\phi_0, \ \pa_t \phi(0)=\phi_1, \quad x\in\R^2.
  \end{array}
 \end{array}\right.
\label{eq_free}
\end{align}
\begin{lem}\label{lem1}
For $\phi_0, \phi_1 \in C_0^{\infty}(\R^2)$ and $\alpha\in\Z_+^3$, 
there is a positive constant $C=C_{\alpha}(\phi_0, \phi_1)$ such that 
the smooth solution $\phi$ to \eqref{eq_free} satisfies
\begin{align}\label{lem_free}
\left| \pa^{\alpha} \phi(t,x) \right| \le C\jb{t+|x|}^{-1}\jb{t-|x|}^{-|\alpha|-1/2}, \qquad 
(t,x)\in [0,\infty)\times \R^2.
\end{align}
\end{lem}
\begin{lem}\label{lem2}
For $\phi_0, \phi_1 \in C_0^{\infty}(\R^2)$, 
there is a positive constant $C=C(\phi_0, \phi_1)$ such that 
the smooth solution $\phi$ to \eqref{eq_free} satisfies
\begin{align}\label{lem_free_rad}
\left| |x|^{1/2}\pa_a \phi(t,x) - \omega_a(\pa_{\sigma}\op{F}_0[\phi_0,\phi_1])(|x|-t,\omega) \right| 
\le C\jb{t+|x|}^{-1}\jb{t-|x|}^{-1/2}
\end{align}
for $(t,x)\in [0,\infty)\times \R^2\setminus\{0\}$, 
$a=0,1,2$, with the convention $\omega_0=-1$, $\omega_1=x_1/|x|$, $\omega_2=x_2/|x|$.
\end{lem}

\begin{lem}\label{lem3}
For $\phi_0, \phi_1 \in C_0^{\infty}(\R^2)$, 
there is a positive constant $C=C(\phi_0, \phi_1)$ such that 
\begin{align}\label{lem_rad}
\left| \pa_{\sigma}\mathcal{F}_0[\phi_0,\phi_1](\sigma,\omega) \right| 
\le C\jb{\sigma}^{-3/2}, \qquad
(\sigma,\omega)\in \R\times \Sph^1.
\end{align}
\end{lem}

For the proof of Lemmas~\ref{lem1}, \ref{lem2} and \ref{lem3}, see Section~3 in \cite{K}.

\section{Proof of the main results}  \label{sec_proofs}
In this section, we prove the Theorem~\ref{thm_main} and Corollary~\ref{cor_1}.
In what follows, we denote several positive constants by the same letter $C$, 
which may be different from one line to another.
\subsection{Proof of Theorem~\ref{thm_main}}  \label{sec_proof_thm}

First we note that $t_{0,\sigma}$ in \eqref{m} can be replaced by 
$t_{1,\sigma}:=\max \{ \eps^{-1},-2\sigma \}$ since we have 
\begin{align}\label{t_1}
\bigl( V_1(t_{1,\sigma})^2 &- V_2(t_{1,\sigma})^2 \bigr)
- \lr{V_1(t_{0,\sigma})^2-V_2(t_{0,\sigma})^2} \notag \\
&=2\int_{t_{0,\sigma}}^{t_{1,\sigma}} 
\lr{V_1(\tau)\pa_t V_1(\tau)-V_2(\tau)\pa_t V_2(\tau)}  \, d\tau \notag \\
&=2\int_{t_{0,\sigma}}^{t_{1,\sigma}} 
\rho(\tau;\sigma,\omega) \, d\tau.
\end{align}
We also recall the following estimate for $\rho$ obtained in Section 4 in \cite{NS}:
\begin{align}\label{rho}
\LR{\int_{t_{1,\sigma}}^{\infty} \rho(\tau;\sigma,\omega) \, d\tau }
&\le C\eps^{2}\jb{\sigma}^{-3/2} \int_{\eps^{-1}}^{\infty} \tau^{2\mu-3/2} \, d\tau \notag\\
&\le C\eps^{5/2-2\mu}.
\end{align}
From \eqref{m}, \eqref{t_1} and \eqref{rho}, we get
\[
\LR{m(\sigma,\omega) 
- \lr{(V_1(t_{1,\sigma};\sigma,\omega))^2-(V_2(t_{1,\sigma};\sigma,\omega))^2}}
\le C\eps^{5/2-2\mu}
\]
for $(\sigma,\omega)\in\R\times\Sph^1$. 
Thus, to prove Theorem~\ref{thm_main}, it suffices to show 
\begin{align}\label{key}
V_j(t_{1,\sigma};\sigma,\omega) 
= \eps \pa_{\sigma}\mathcal{F}_j(\sigma,\omega) + O(\eps^{2-\mu})
\end{align}
as $\eps\to+0$ uniformly in $(\sigma,\omega)\in\R\times\Sph^1$ for $j=1,2$.
The rest part of this subsection is devoted to the proof of \eqref{key}.
We divide the argument into the following two cases.
\begin{itemize}
\item
\noindent{\bf \underline{Case 1: $\sigma\le -1/(2\eps)$.}}\ \ 
If we assume $|x|\le t/2$ and $t\ge \eps^{-1}$, 
we have $\eps^{-1}\le t\le\jb{t+|x|}\le C\jb{t-|x|}$.
It follows from \eqref{apriori_u} and \eqref{apriori_pa_u} that
\begin{align*}
\LR{U(t,x)} \notag
&\le C|x|^{-1/2}\LR{u(t,x)} +C|x|^{1/2}\LR{\pa u(t,x)} \\ \notag
&\le C\eps |x|^{-1/2}\jb{t+|x|}^{-1/2+\mu} 
+ C\eps |x|^{1/2}\jb{t+|x|}^{-1/2}\jb{t-|x|}^{\mu-1}\\
&\le C\eps^{2-\mu}
\end{align*}
for $|x|\le t/2$ and $t\ge \eps^{-1}$. 
Then we obtain
\begin{align}
\LR{V(t;\sigma,\omega)}
=|U(t,(t+\sigma)\omega)|
\le C\eps^{2-\mu} \label{V_away}
\end{align}
for $t+\sigma\le t/2$ and $t\ge \eps^{-1}$. 
In the case $\eps^{-1}\le -2\sigma$, 
we have $t_{1,\sigma}+\sigma=t_{1,\sigma}/2$ and 
$t_{1,\sigma}\ge\eps^{-1}$.
Therefore, from \eqref{lem_rad}, \eqref{V_away} 
and $|\sigma|\ge 1/(2\eps)$, 
we get
\begin{align*}
|V_j(t_{1,\sigma};\sigma,\omega) 
- \eps \pa_{\sigma}\op{F}_j(\sigma,\omega) | \notag
&\le |V_j(t_{1,\sigma};\sigma,\omega)| 
+ \eps |\pa_{\sigma}\op{F}_j(\sigma,\omega) | \\ \notag
&\le C\eps^{2-\mu} + C\eps \jb{\sigma}^{-3/2}\\
&\le C\eps^{2-\mu}.
\end{align*}

\item

\noindent{\bf \underline{Case 2: $\sigma> -1/(2\eps)$.}}\ \ 
For $j=1,2$, let $u^0_j=u^0_j(t,x)$ be  the solution to the free wave equation $\Box u^0_j = 0$ 
with the initial data $u^0_j(0)=f_j$, $\pa_t u^0_j(0)=g_j$ 
and we put $u^1_j(t,x):=u_j(t,x)-\eps u^0_j(t,x)$, so that $u_j^1$ solves
\begin{align*}
&\Box u^1_j(t,x) = F_j(\pa u), \ \quad\quad\quad (t,x)\in (0,\infty)\times\R^2, \\
&u^1_j(0,x)=\pa u^1_j(0,x)=0, \quad x\in\R^2.
\end{align*}
We also define $U^{l}(t,x)=\mathcal{D}(|x|^{1/2}u^{l}(t,x))$ and 
$V^{l}(t;\sigma,\omega):=U^{l}(t,(t+\sigma)\omega)$, for $l=0,1$, respectively.
It follows from \eqref{lem_free} and \eqref{lem_free_rad} that
\begin{align*}
&\LR{U_j^0(t,x) - \pa_{\sigma}\op{F}_j(|x|-t,x/|x|)} \notag\\
\le &\frac{1}{2}\sum_{a=0}^{2}||x|^{1/2}\pa_a u_j^0(t,x) 
- \omega_a \pa_{\sigma}\op{F}_j(|x|-t,\omega)| 
+ \frac{1}{4|x|^{1/2}}| u_j^0(t,x) |\notag\\
\le &C\jb{t+|x|}^{-1}\jb{t-|x|}^{-1/2} 
+ C|x|^{-1/2}\jb{t+|x|}^{-1/2}\jb{t-|x|}^{-1/2} \notag \\
\le &C\eps
\end{align*}
for $|x|\ge1/(2\eps)$. 
Hence we get
\begin{align}\label{V^0}
\LR{V_j^0(t;\sigma,\omega) - \pa_{\sigma}\op{F}_j(\sigma,\omega)}
\le C\eps
\end{align}
for $t+\sigma\ge1/(2\eps)$.
We next consider the estimate for $V^1$. 
Note that we have 
$(\Gamma^{\alpha} \phi,\pa_t\Gamma^{\alpha} \phi)|_{t=0}
\in \lr{C_0^{\infty}(\R^2)}^2$ and 
$\|\Gamma^{\alpha} \phi(0)\|_{L^{\infty}(\R^2)}$, 
$\|\pa_t\Gamma^{\alpha} \phi(0)\|_{L^{\infty}(\R^2)}=O(\eps^3)$
for $\alpha\in\Z_+^7$ if $\phi(t,x)$ satisfies $\Box \phi = N(\pa\phi)$ 
with a cubic nonlinear term $N(\pa\phi)$ and
$(\phi,\pa_t\phi)|_{t=0}\in \lr{C_0^{\infty}(\R^2)}^2$.
By using \eqref{apriori_pa_u_k}, \eqref{apriori_energy} and the standard energy method 
for $\Gamma^{\alpha} u^1$ with $|\alpha|\le2$, we obtain
\begin{align*}
\| \pa u^{1}(t) \|_{2} 
&\le C\eps^3 
+ C\int_0^t \left\| |\pa u(\tau,\cdot)|_1\right\|_{L^{\infty}}^2 \| \pa u(\tau) \|_2 \, d\tau \notag \\
&\le C\eps^3 + C\eps^3 \int_0^{\eps^{-1}} (1+\tau)^{-1+\mu+\nu} \, d\tau \notag \\
&\le C\eps^3 + C\eps^3(1+\eps^{-1})^{\mu+\nu} \notag \\
&\le C\eps^{3-\mu-\nu}
\end{align*}
for $0\le t \le \eps^{-1}$.
Then, by the Klainerman-Sobolev inequality, we get
\begin{align}\label{1}
\jb{t+|x|}^{1/2}|\pa u^1(t,x)|\le C\eps^{3-\mu-\nu}
\end{align}
for $0\le t \le \eps^{-1}$, $x\in \R^2$.
It follows from \eqref{apriori_u} and \eqref{lem_free} that
\begin{align}\label{0}
|x|^{-1/2}\LR{u^{1}(t,x)}
&\le |x|^{-1/2}\lr{|u(t,x)| + \eps|u^{0}(t,x)|} \notag \\
&\le |x|^{-1/2}\lr{C\eps \jb{t+|x|}^{\mu-1/2} 
+ C\eps\jb{t+|x|}^{-1/2}\jb{t-|x|}^{-1/2}} \notag \\
&\le C\eps^{2-\mu}
\end{align}
for $|x|\ge 1/(2\eps)$.
From \eqref{1} and \eqref{0}, we get 
\begin{align*}
\LR{U^{1}(t,x)}
&\le C|x|^{1/2}|\pa u^{1}(t,x)| + C|x|^{-1/2}|u^{1}(t,x)| \notag \\
&\le C\eps^{3-\mu-\nu} + C\eps^{2-\mu} \notag \\
&\le C\eps^{2-\mu}
\end{align*}
for $|x|\ge 1/(2\eps)$, $0\le t\le \eps^{-1}$.
Therefore, we obtain 
\begin{align}\label{V^1}
\LR{V^{1}(t;\sigma,\omega)}
\le C\eps^{2-\mu}
\end{align}
for $t+\sigma\ge 1/(2\eps)$, $0\le t\le \eps^{-1}$.
When $\eps^{-1}> -2\sigma$, we have $t_{1,\sigma} = \eps^{-1}$ 
and $t_{1,\sigma}+\sigma>t_{1,\sigma}/2= 1/(2\eps)$. 
Thus, by \eqref{V^0} and \eqref{V^1}, we get
\begin{align*}
|V_j(t_{1,\sigma};\sigma,\omega) - \eps \pa_{\sigma}\op{F}_j(\sigma,\omega) | \notag
&\le |V_j^1(t_{1,\sigma};\sigma,\omega)| 
+ \eps |V_j^0(t_{1,\sigma};\sigma,\omega) - \pa_{\sigma}\op{F}_j(\sigma,\omega) | \\ \notag
&\le C\eps^{2-\mu}.
\end{align*}
\end{itemize}
Combining the two cases above, we arrive at the desired expression \eqref{key}. 
This completes the proof of Theorem~\ref{thm_main}.
\qed
\subsection{Proof of Corollary~\ref{cor_1}}  \label{sec_proof_cor_1}
Let $E=\{ (\sigma,\omega)\in\R\times\Sph^1 ; 
|\pa_{\sigma}\op{F}_1(\sigma,\omega)|>|\pa_{\sigma}\op{F}_2(\sigma,\omega)| \}$.
By \eqref{1>2}, $E$ is a non-empty open set.
We can take a bounded open set $\op{M}$ in $\R$ and 
an open set $\op{N}$ in $\Sph^1$ such that $\sigma^*\in\op{M}$, $\omega^*\in\op{N}$ 
and $\overline{\op{M}\times\op{N}}\subset E$, where $\overline{\op{M}\times\op{N}}$ 
denotes the closure of $\op{M}\times\op{N}$ in $\R\times\Sph^1$.
Now we put $F=\overline{\op{M}\times\op{N}}$ and 
\[
C_1=\min_{(\sigma,\omega)\in F}
\lr{
\lr{\pa_{\sigma}\op{F}_1(\sigma,\omega)}^2
-\lr{\pa_{\sigma}\op{F}_2(\sigma,\omega)}^2
}.
\]
Then we see that $F$ is compact, and thus $C_1>0$. 
By Theorem~\ref{thm_main}, we have
\[
m(\sigma,\omega)\ge C_1\eps^2-C\eps^{5/2-2\mu}>0
\]
for $(\sigma,\omega)\in F$, if $\eps>0$ is small enough. 
This and Proposition~\ref{prop_NS3} imply
$V_1^+(\sigma,\omega)\neq 0$ for $(\sigma,\omega)\in F$, 
whence $\| V_1^+ \|_{L^2(F)}>0$. By virtue of \eqref{NS2}, 
we can take $T_1>0$ such that 
\[
\| \pa u_1(t,\cdot) - \hat\omega(\cdot)V_1^{+,\#}(t,\cdot) \|_{L^2(\R^2)}
< \frac{1}{\sqrt{2}} \| V_1^+ \|_{L^2(F)}
\]
for $t>T_1$. Therefore we have
\begin{align*}
\| u_1(t) \|_E
&\ge \lr{\frac{1}{2}\int_{\R^2} |\hat\omega(x)V_1^{+,\#}(t,x) |^2\, dx }^{1/2}
- \lr{\frac{1}{2}\int_{\R^2} | \pa u_1(t,x) - \hat\omega(x)V_1^{+,\#}(t,x) |^2\, dx}^{1/2} \\
&= \| V_1^{+,\#}(t,\cdot) \|_{L^2(\R^2)}- \frac{1}{\sqrt{2}}\| 
\pa u_1(t,\cdot) - \hat\omega(\cdot)V_1^{+,\#}(t,\cdot) \|_{L^2(\R^2)} \\
&\ge \| V_1^{+} \|_{L^2(F)}- \frac{1}{2}\| V_1^{+} \|_{L^2(F)} \\
&=\frac{1}{2}\| V_1^{+} \|_{L^2(F)}
\end{align*}
for $t>T_1$. Consequently, 
\[
\lim_{t\to+\infty}\| u_1(t) \|_E \ge  \frac{1}{2}\| V_1^{+} \|_{L^2(F)} >0,
\]
as desired. 
Interchanging the role of $u_1$ and $u_2$, we also have 
$\lim_{t\to+\infty}\| u_2(t) \|_E >0$.
\qed

\medskip
\subsection*{Acknowledgments}
The author would like to thank Professor Soichiro Katayama and 
Professor Hideaki Sunagawa for their useful conversations on this subject.
This work was partly supported by Osaka City University Advanced Mathematical 
Institute (MEXT Joint Usage/Research Center on Mathematics and Theoretical 
Physics JPMXP0619217849).


\end{document}